\documentclass[12pt]{amsart}
\usepackage{amsmath}
\usepackage{amssymb}
\usepackage[mathcal]{eucal}
\usepackage[all]{xy}
\newtheorem{deff}{Definition}[section]

\newtheorem{theorem}[deff]{Theorem}
\newtheorem{corollary}[deff]{Corollary}

\newtheorem{proposition}[deff]{Proposition}

\newtheorem{observation}[deff]{Observation}

\newtheorem{em-example}[deff]{Example}
\newtheorem{em-def}[deff]{Definition}        
\newtheorem{em-remark}[deff]{Remark}         
\newtheorem{em-question}[deff]{Question}

\newtheorem{problem}[deff]{Problem}

\newenvironment{example}{\begin{em-example} \em }{ \end{em-example}}

\newenvironment{remark}{\begin{em-remark} \em }{\end{em-remark}}

\newcommand{\lcal}{\mathcal {L}}

\newcommand{\ucal}{\mathcal {U}}
\newcommand{\vcal}{\mathcal {V}}

\def\ker{\mathop{\rm ker}}

\def\cl{\mathop{\it cl}}

\DeclareMathSymbol{\res}{\mathord}{AMSa}{"16}

\def\:{\nobreak \hskip .1111em\mathpunct {}\nonscript \mkern
   -\thinmuskip {:}\hskip .3333emplus.0555em\relax}

\catcode`\@=12

\def\N{{\mathbb N}}
\def\R{{\mathbb R}}

\title[On linear continuous operators between distinguished spaces $C_p(X)$]{On linear continuous operators between distinguished spaces $C_p(X)$}
\author{Jerzy K\c{a}kol and Arkady Leiderman}

\address{Faculty of Mathematics and Informatics, A. Mickiewicz University,
61-614 Pozna\'{n}, Poland and Institute of Mathematics Czech Academy of Sciences, Prague, Czech Republic}
\email{kakol@amu.edu.pl}

\address{Department of Mathematics, Ben-Gurion University of the Negev, Beer Sheva, P.O.B. 653, Israel}
\email{arkady@math.bgu.ac.il}

\keywords{Distinguished locally convex space, linear continuous operator, $\Delta$-space, countable ordinal}
\subjclass[2010]{Primary 54C35; Secondary 46A03, 46A20}
\date{\today}

\begin{document}

\begin{abstract}
As proved in \cite{ka-arkady}, for a Tychonoff space $X$, a locally convex space $C_{p}(X)$ is distinguished if and only if $X$ is a $\Delta$-space.
If there exists a linear continuous surjective mapping $T:C_p(X) \to C_p(Y)$ and $C_p(X)$ is distinguished,
 then $C_p(Y)$ also is distinguished \cite{Kakol-Leiderman}.

Firstly, in this paper we explore the following question: Under which conditions the operator $T:C_p(X) \to C_p(Y)$ above is open?
Secondly, we devote a special attention to concrete distinguished spaces $C_p([1,\alpha])$, where $\alpha$ is a countable ordinal number.
A complete characterization of all $Y$ which admit a linear continuous surjective mapping $T:C_p([1,\alpha]) \to C_p(Y)$ is given.
We also observe that for every countable ordinal $\alpha$ all closed linear subspaces of
$C_p([1,\alpha])$ are distinguished, thereby answering an open question posed in \cite{Kakol-Leiderman}.

Using some properties of $\Delta$-spaces we prove that
a linear continuous surjection  $T:C_p(X) \to C_k(X)_w$, where $C_k(X)_w$ denotes the Banach space $C(X)$ endowed with its weak topology,
does not exist
for every infinite metrizable compact $C$-space $X$ (in particular, for every infinite compact $X \subset \R^n$).
\end{abstract}

\thanks{The first named author
 is supported by the GA\v{C}R project 20-22230L and RVO: 67985840.\\
The authors thank J.C. Ferrando for providing a short argument in
Remark \ref{rem1} (iii) and W.B. Johnson and T. Kania for a very helpful advice and discussion about Theorem \ref{cpp}.}
\maketitle

\section{Introduction}\label{intro}

A locally convex space (lcs) $E$  is called \emph{distinguished } if  its strong dual $E^{\prime}_{\beta}=(E^{\prime},\beta(E^{\prime},E))$ is a barrelled space.
A. Grothendieck \cite{Gro}  proved  that  a metrizable lcs is distinguished if and only if $E^{\prime}_{\beta}$ is bornological.
Also, if all bounded subsets of the strong dual $E^{\prime}_{\beta}$ of a metrizable lcs are metrizable,
 then $E$ is distinguished \cite{Gro}.
Recall that the \emph{strong topology} $\beta$ on $E^{\prime}$ is the topology of uniform convergence on bounded subsets of $E$.
A subset $A\subset E$ is \emph{bounded} if is absorbed by each neighbourhood of zero in $E$.
 We refer the reader to the survey article  of K. D. Bierstedt and J. Bonet \cite{BB2}, where
many classes of distinguished Fr\'echet lcs are presented (reflexive spaces, Montel spaces, nuclear spaces, Schwartz spaces, etc).

In this paper we continue the study of distinguished lcs in the frame of spaces $C_p(X)$ and $C_k(X)$ (developed earlier in \cite{FK}, \cite{FKLS}, \cite{FS},
 \cite{ka-arkady}, \cite{Kakol-Leiderman}).
By $C_{p}(X)$ and $C_{k}(X)$  we mean the spaces of all real-valued continuous functions on a Tychonoff space $X$
 endowed with the topology of pointwise convergence $\tau_{p}$ and the compact-open topology $\tau_{k}$, respectively.

Recall that for a vector space $E$ the finest locally convex topology $\xi$ of $E$ is generated by the family of  all absolutely convex and absorbing subsets of $E$
which form a base of neighbourhoods of zero for the topology $\xi$. We denote by $\N$ the infinite countable discrete space.
It is obvious that  $C_{p}(\N)=\mathbb{R}^{\N}$ is the only distinguished  Fr\'echet lcs $E$ for which its strong dual carries  the finest locally convex topology.
 In fact, $C_{p}(\N)^{\prime}_{\beta}$ is isomorphic to $\varphi$, i.e. the $\aleph_{0}$-dimensional vector space with the finest locally convex topology.
On the other hand, one can characterize distinguished lcs $C_p(X)$ as follows: $C_{p}(X)$ is distinguished if and only if $C_{p}(X)^{\prime}_{\beta}$ carries the finest locally convex topology,
 see \cite{FK}, \cite{FKLS}.

We proved that $C_{p}(X)$ is distinguished if and only if $X$ is a $\Delta$-space \cite{ka-arkady}.\footnote{We should mention
that independently and simultaneously an analogous description of distinguished $C_p$-spaces (but formulated in different terms)
appeared in \cite{FS}.}
The class $\Delta$ of $\Delta$-spaces naturally extends the class of $\Delta$-sets of reals.
Note that the original definition of a $\Delta$-set $A \subset \R$ is due to G. Reed and E. van Douwen \cite{Reed}.

A topological space $X$ is said to be a \emph{$\Delta$-space} if for every decreasing sequence $\{D_n: n \in \omega\}$
of subsets of $X$ with empty intersection, there is a decreasing sequence $\{V_n: n \in \omega\}$ consisting of open
subsets of $X$, also with empty intersection, and such that $D_n \subset V_n$ for every $n \in \omega$ \cite{ka-arkady}.

Quite recently a range of new facts about $\Delta$-spaces has been obtained by the authors.
We proved that every \v{C}ech-complete (in particular, compact) $X\in\Delta$ is scattered \cite{ka-arkady}.
Moreover, every countably compact space is scattered \cite{Kakol-Leiderman}.
Class $\Delta$ contains all separable compact spaces of the Isbell-Mr\'owka type, while there are compact scattered spaces $X\notin\Delta$, for example, ordinal space $[0, \omega_{1}]$ \cite{ka-arkady}.
In fact, every compact $\Delta$-space must have countable tightness \cite{Kakol-Leiderman}.
Although the property of being  a $\Delta$-space is not inherited by continuous images,
the class $\Delta$ does preserve closed continuous images \cite{ka-arkady}.
Furthermore, a countable union of compact $\Delta$-spaces is also a $\Delta$-space;
 in particular, $\sigma$-product of any family consisting of scattered Eberlein compact spaces is a $\Delta$-space \cite{Kakol-Leiderman}.
The next result is a matter of high importance to us.

\begin{theorem}\cite{Kakol-Leiderman}\label{th:dom}
The class $\Delta$ is invariant under the relation of $l$-dominance, i.e. if $C_p(X)$ is distinguished
and there exists a linear continuous surjective mapping $T:C_p(X) \to C_p(Y)$, then $C_p(Y)$ also is distinguished.
\end{theorem}

The focus in this paper is on the following question: Under which conditions the operator $T:C_p(X) \to C_p(Y)$ above is open?
We notice that in general there are plenty of linear continuous non-open surjections between $C_p(X)$-spaces.

\begin{proposition}\cite{Leiderman-1}\label{prop:notopen}
Let $X$ be a Tychonoff space and $h: X \to X$ be
a continuous mapping such that the orbit of at least one point
is infinite. Then the mapping
 $T: C_p(X) \to C_p(X)$ defined by
$Tf(x) = \lambda f(x) + f(h(x))$ is a linear continuous non-open surjection for each $\lambda$ with
$\vert\lambda\vert > 1$.
\end{proposition}

In Section \ref{sec2} we prove  the following open mapping theorem for spaces $C_{p}(X)$ over $\Delta$-spaces $X$.
\begin{theorem}\label{some-remark}
Let $X$ be a $\Delta$-space.
Consider the following assertions for a linear continuous surjection $T: C_p(X)\rightarrow C_{p}(Y)$.
\begin{enumerate}
\item [(1)] $C_{p}(X)$ admits  a  base $\mathcal{U}$ of absolutely convex  neighbourhoods of zero such that $\cl{T(U)}\subset 2T(U)$ for each $U\in\mathcal{U}$.
\item [(2)] $T$ is an open mapping.
\item [(3)] $C_p(X)$ admits a base $\mathcal{U}$ of absolutely convex closed neighbourhoods of zero such that  $T(U)$ is closed  for each $U\in\mathcal{U}$.
\end{enumerate}
Then (3) $\Rightarrow$ (2) $\Leftrightarrow$  (1).
\end{theorem}

\begin{corollary}\label{cor1}
If $X$ and $Y$ are $\Delta$-spaces and $T:C_{p}(X)\rightarrow C_{p}(Y)$ is a linear continuous bijection, then $T$ is an isomorphism if and only if any condition from Theorem \ref{some-remark} holds.
\end{corollary}

It is easy to see that in above Proposition \ref{prop:notopen}, $X$ might be selected as a $\Delta$-space for which a continuous linear non-open surjection
$T:C_{p}(X)\to C_{p}(X)$ exists but $T$ lacks condition (3) of Theorem \ref{some-remark}. This happens, for example, if $C_{p}(X)$ is barrelled.
Clearly, there exist Banach spaces $E$ and $F$ and  linear continuous surjections $T:E\rightarrow F$ with non-closed  $T(S)$ for the closed unit ball $S\subset E$;
 injective linear continuous maps $T:E\rightarrow F$ such that $T(S)$ is closed  were studied under the name  \emph{semi-embedding}, see \cite{Peck}.

We stress that our proof of Theorem \ref{some-remark} rely inevitably on the assumption that $X$ is a $\Delta$-space.
Recall that the unit segment $[0,1]$ is not a $\Delta$-space.
A challenging problem of whether
a non-open continuous linear surjection $T:C_p[0,1]\rightarrow C_p[0,1]$ satisfying condition (3) of Theorem \ref{some-remark} exists,
remains unsolved.

We devote a special attention to lcs $C_p(X)$ for some concrete $\Delta$-spaces $X$. For any ordinal number $\alpha$ let $[1,\alpha]$ denote
the ordered compact space consisting of all ordinals $\gamma \leq \alpha$.
The symbols $\omega$ and $\omega_1$ stand for the first infinite and the first uncountable ordinal numbers, respectively.
S. Gul'ko showed that for any pair of infinite countable ordinals $\alpha$ and $\beta$, the spaces $C_{p}([1,\alpha])$ and $C_{p}([1,\beta])$ are uniformly homeomorphic, see \cite{Gul'ko1}.
 Theorem \ref{cpp} shows however that if $\alpha \leq \beta < \omega_1$ and
the spaces $C_{p}([1,\alpha])$ and $C_{p}([1,\beta])$ are not linearly isomorphic, then $C_{p}([1,\beta])$ is not a continuous linear image of $C_{p}([1,\alpha])$.
Thus our results are related to the research area of linear topological classification of function spaces
 originated long time ago in the pioneering work of C. Bessaga and A. Pe{\l}czy\'nski \cite{bessaga}.
A complete linear topological classification of the spaces $C_p([1,\alpha])$ has been obtained independently in
\cite{Gul'ko} (for all $\alpha$) and \cite{Baars} (for countable $\alpha$).

We also observe that for every countable ordinal compact space $X=[1,\alpha]$ the space $C_p(X)$ is hereditarily distinguished in the following sense:
all closed linear subspaces of $C_p(X)$ are distinguished (Theorem \ref{th:hered}). This result answers in the affirmative Problem 2.18 posed in \cite{Kakol-Leiderman}.

The argument used in the proof of Theorem \ref{some-remark} applies to provide a "linear version" of  M. Krupski and W. Marciszewski results from \cite{Krupski-1}, \cite{Krupski-2}
 about the non-existence of a homeomorphism between spaces $C_{p}(X)$ and $C_k(Y)_{w}$ over certain compact spaces $X$ and $Y$,
where $C_k(Y)_{w}$ denotes the Banach space $C(Y)$ endowed with its weak topology $w$.
Although we do not know if there exists an infinite compact space $X$ admitting a linear continuous surjection $T:C_{p}(X)\rightarrow C_k(X)_{w}$, we show
 (Theorem \ref{the-1}) that any infinite metrizable compact $C$-space $X$, (in particular, any infinite metrizable finite-dimensional compact space $X$),
 does not admit such a mapping.

Our notations are standard, the reader is advised to consult with the monograph
\cite{Jarchow} and the survey paper \cite{Rosenthal} for the notions which are not explicitly defined in the text.
In the article we pose several open questions.

\section{Proof of Theorem \ref{some-remark} and  $\varphi$-dual subspaces}\label{sec2}

A lcs $E$ is called \emph{free} if $E$ carries the finest locally convex topology.
Each free lcs has only finite-dimensional bounded sets, see for example \cite[Theorem 2.1]{ka-ta-joh}. On the other hand, for every Tychonoff $X$,
in the space $C_{p}(X)'_{\beta}$ all bounded sets are finite-dimensional, but $C_{p}(X)'_{\beta}$ is not necessarily free.
For the proofs of the next characterizations we refer to \cite{FKLS} and \cite{ka-arkady}.

\begin{theorem} \label{the1}  Let $X$ be a Tychonoff space.
\begin{enumerate}
\item [(a)] \cite{FK}, \cite{FKLS} $C_{p}(X)$ is distinguished if and only if $C_{p}(X)'_{\beta}$ is free.
\item [(b)] \cite{FKLS} $C_{p}(X)$ is distinguished if and only if for each  bounded set $A\subset\mathbb{R}^{X}$
 there exists a bounded set $B\subset C_{p}(X)$ such that $A \subset \cl_{\mathbb{R}^{X}}(B)$.
\item [(c)] \cite{ka-arkady} $C_{p}(X)$ is distinguished if and only if $X$ is a $\Delta$-space.
\end{enumerate}
\end{theorem}

We start with the following theorem describing the free spaces $C_{k}(X)'_{\beta}$ which (together with Corollary \ref{Ck-1}) will be used in proving
 Theorem \ref{some-remark}.
By $C_{p}(X)'$ and $C_{k}(X)'$ we denote the topological dual of $C_{p}(X)$ and $C_{k}(X)$, respectively.

\begin{theorem}\label{Ck}
Let $X$ be a Tychonoff space. Let $\vartheta_{X}$ be the strong topology of  $C_{k}(X)'_{\beta}$. Then the following assertions are equivalent:
\begin{enumerate}
\item [(1)] The topology $\vartheta_{X}$  is the finest locally convex topology of $C_{k}(X)'_{\beta}$.
\item [(2)] $\tau_{k}=\tau_{p}$ and $C_{p}(X)$ is distinguished.
\item [(3)] For each bounded set $A\subset\mathbb{R}^{X}$ there exists a bounded set $B$ in $C_{k}(X)$ such that $A \subset \cl_{\mathbb{R}^{X}}(B)$.
\end{enumerate}
\end{theorem}

\begin{proof} Implications (2) $\Rightarrow$ (1) and (2) $\Rightarrow$ (3) are straightforward consequences of Theorem \ref{the1}.
(1) $\Rightarrow$ (2):  Note  that the topology $\tau_{k}$ on $C(X)$ is $\tau_{p}$-polar, i.e. $\tau_{k}$ admits a base of absolutely convex neighbourhoods of zero consisting of $\tau_{p}$-closed sets.  Hence, if $S$ is a $\tau_{k}$-bounded set  in $C_{k}(X)$, then
 the closure $\cl_{\tau_{p}}(S)$ is still $\tau_{k}$-bounded.

 Assume that (1) holds. Since $\tau_{p}\leq\tau_{k}$, the strong  topology $\beta_{X}$ of $C_{p}(X)'$ is finer than the strong topology
$\vartheta_{X}$ of $C_{k}(X)'$  restricted to  $C_{p}(X)'$ (shortly $\vartheta_{X}|C_{p}(X)'\leq\beta_{X}$). As the topology $\vartheta_{X}$ is the finest locally convex topology on $C_{k}(X)'$, we have in fact the equality  $\vartheta_{X}|C_{p}(X)'=\beta_{X}$. Hence $(C_{p}(X)',\beta_{X})$ is a topological subspace of
 $(C_{k}(X)',\vartheta_{X})$ and carries the finest locally convex topology, so $C_{p}(X)$ is distinguished by Theorem \ref{the1}.
 Next we show that $\tau_{p}=\tau_{k}$.

We know that for every absolutely convex neighbourhood of zero $U$ in the space  $(C_{p}(X)',\beta_{X})$ there exists an absolutely convex neighbourhood of
zero $V$ in $(C_{k}(X)',\vartheta_{X})$ such that $V\cap C_{p}(X)'\subset U$.
 Consequently, for each $\tau_{p}$-closed bounded absolutely convex set $B\subset C_{p}(X)$ there exists an  absolutely convex $\tau_{k}$-bounded set
 $S\subset C_{k}(X)$ with $S^{\bullet}\cap C_{p}(X)'\subset B^{\circ}$, where
$$S^{\bullet}=\{x^*\in C_{k}(X)': |x^*(f)|\leq 1, f\in S\},$$
$$B^{\circ}=\{x^*\in C_{p}(X)':|x^*(f)|\leq 1, f\in B\}.$$
 Hence, $\{x^*\in C_{p}(X)':|x^*(f)|\leq 1, f\in \cl_{\tau_{p}}(S)\}\subset
\{x^*\in C_{p}(X)':|x^*(f)|\leq 1, f\in S\}\subset \{x^*\in C_{p}(X)':|x^*(f)|\leq 1, f\in B\},$
what shows that  the polar set $B^{\circ\circ}= \{x^*\in C_{p}(X)':|x^*(f)|\leq 1, f\in B\}^{\circ}$ is contained in  the polar set
$(\cl_{\tau_{p}}(S))^{\circ\circ}= \{x^*\in C_{p}(X)':|x^*(f)|\leq 1, f\in \cl_{\tau_{p}}(S)\}^{\circ}$, where the corresponding polars are taken in the
 dual pair $(C_{p}(X)', C_{p}(X))$.

By the bipolar theorem \cite[Theorem 8.2.2]{Jarchow}, $B = B^{\circ\circ}\subset (\cl_{\tau_{p}}(S))^{\circ\circ}=\cl_{\tau_{p}}(S)$.
  Since $\cl_{\tau_{p}}(S)$ is $\tau_{k}$-bounded,
 every $\tau_{p}$-bounded set in $C(X)$ is $\tau_{k}$-bounded.
On the other hand,  $C_{p}(X)$ is \emph{quasibarrelled}, i.e. every absolutely convex closed set absorbing bounded sets in $C_{p}(X)$ is a neighbourhood of zero, see \cite[11.7.3]{Jarchow}.
 Now take a $\tau_{p}$-closed  absolutely convex neighbourhood of zero $U$  in $\tau_{k}$. Since, as we proved,  every $\tau_{p}$-bounded set is $\tau_{k}$-bounded, the set $U$ absorbs bounded sets
in $C_{p}(X)$. Hence, $U$ is also a neighbourhood of zero in $\tau_{p}$, so $\tau_{p}=\tau_{k}$.

(3) $\Rightarrow$ (2): Take a bounded set $A\subset C_{p}(X)$. Then there exists a bounded set $B\subset C_{k}(X)$ with
 $A \subset \cl_{\mathbb{R}^{X}}(B) \cap C_{p}(X)(\subset \cl_{\tau_{p}}(B)$).
 Clearly, the closure $\cl_{\tau_{p}}(B)$ is still $\tau_{k}$-bounded.  Therefore,  the topologies $\tau_{p}$ and $\tau_{k}$ have the same bounded sets,
 and then $\tau_{p}=\tau_{k}$. Applying Theorem \ref{the1} we deduce that $C_{p}(X)$ is distinguished.
\end{proof}

A similar argument, as used in the  proof of Theorem \ref{Ck} ((1)$\Rightarrow$ (2))  for the case  $\xi=\tau_{k}$, applies for the next Corollary \ref{Ck-1}.
 Analogously, one shows that bounded sets in both topologies $\xi$ and $\tau_{p}$ are the same. Since $C_{p}(X)$ is quasibarrelled,  $\tau_{p}$-polarity of $\xi$
 yields the equality  $\xi=\tau_{p}$.

\begin{corollary}\label{Ck-1}
Let  $\xi\geq\tau_{p}$ be a locally convex topology on $C(X)$ which is $\tau_{p}$-polar, i.e. $\xi$ has a base of absolutely convex $\tau_{p}$-closed  neighbourhoods of zero.
 If $(C(X),\xi)'_{\beta}$ is free, then $\xi=\tau_{p}$.
\end{corollary}

\begin{example}\label{dis}
If $X$ is any uncountable compact space, then the Banach space $C_{k}(X)$ is evidently distinguished but $C_{k}(X)'_{\beta}$ is not free.
\end{example}

\begin{proof} [Proof of Theorem \ref{some-remark}]
(3) $\Rightarrow$  (2): Let $Q:C_{p}(X)\rightarrow C_{p}(X)/\ker(T)$ be the quotient map and  let $\hat{T}:C_{p}(X)/\ker(T)\rightarrow C_{p}(Y)$
be a map associated with $T$. It is well known that  the adjoint map
$Q^{*}: (C_{p}(X)/\ker(T))^{'}_{\beta} \rightarrow C_{p}(X)^{'}_{\beta}$ is continuous and injective.  Since $C_{p}(X)'_{\beta}$ is free,
 the space $(C_{p}(X)/\ker(T))'_{\beta}$ is free, too. Moreover, $\widehat{T}$ is injective and  $\widehat{T}(Q(U))=\widehat{T}(\cl(Q(U)))=T(U)$
is closed in $C_{p}(Y)$ for $U\in\mathcal{U}$.
Hence there exists a locally convex topology $\xi$ on $C(Y)$  which is $\tau_{p}$-polar with  $\tau_{p}\leq\xi$  and
$\widehat{T}: C_{p}(X)/\ker(T)\rightarrow  (C(Y),\xi)$ is an isomorphism. Since $(C(Y),\xi)^{'}_{\beta}$ is free,
 Corollary \ref{Ck-1} applies to get that $\tau_{p}=\xi$, so $T$ is open.

(2) $\Rightarrow$  (1) is clear.  (1) $\Rightarrow$  (2): First note that to get the conclusion of  Corollary \ref{Ck-1} it is enough to assume
that $\xi$ has a base $\mathcal{U}$ of absolutely convex neighbourhoods of zero with
$\cl_{\tau_{p}}(U) \subset 2U$ for each $U\in\mathcal{U}$. This assumption applies to show that $\cl_{\tau_{p}}(B)$ is $\xi$-bounded for each
 $\xi$-bounded $B\subset C(X)$, which is  the  essential  part used in the proof of  Theorem \ref{Ck}.
Next, we argue as in the proof (3)$\Rightarrow$ (2) above.
\end{proof}

\begin{remark}\label{rem1}
(i)  Theorem \ref{some-remark} fails if $C_{p}(X)$ is replaced by $C_{k}(X)$. Indeed, for every infinite compact space $X$, if $T: C_{k}(X)\rightarrow C_{p}(X)$ is the identity map, then $T$ satisfies condition (3) but $T$ is not open.

(ii) If $C_{p}(Y)$ is a \emph{barrelled} lcs, i.e. every closed absolutely convex absorbing set in $C_{p}(Y)$ is a neighbourhood of zero,
(see \cite[Corollary 11.7.6]{Jarchow}),
 then the implication (3) $\Rightarrow$ (2) of Theorem \ref{some-remark} holds for each $C_{p}(X)$.

(iii) The most natural situation when a linear continuous surjective mapping $T: C_p(X)\rightarrow C_p(Y)$ both is open and satisfies condition (3) of Theorem \ref{some-remark}
is  the following: $X$ is normal, $Y$ is a closed subspace of $X$ and $T:f\mapsto f|_{Y}$ is the restriction map.
 We do not need to assume that $C_p(X)$ is distinguished.
For any finite set $A \subset X$ and a real number $r > 0$ denote $U_{A,r}=\{f\in C(X):|f(x)|\leq r\,\, \mbox{for every}\, x\in A\}$. Define the family of sets
 $\mathcal{U}$ consisting of all sets $U_{A,r}$ of such a form. Then $\mathcal{U}$ is a base of absolutely convex closed neighbourhoods of zero in $C_p(X)$.
We claim that $T(U_{A,r})$ is closed in $C_{p}(Y)$ for every $U_{A,r}\in\mathcal{U}$. Indeed, let $B = A \cap Y$, then
applying the Tietze--Urysohn lemma, it is easy to see that
$T(U_{A,r})= \{f\in C(Y):|f(x)|\leq r\,\, \mbox{for every}\, x\in B\}$.
It follows that $T$ satisfies condition (3).
\end{remark}

We complete this section with a theorem illustrating a role of  distinguished spaces $C_p([1,\omega])$ and $C_p(\N)$.
In what follows, if lcs $E_1$ are $E_2$ are \emph{isomorphic}, we write  $E_1\approx E_2$.
A vector subspace $F$ of $E$ will be called \emph{$\varphi$-dual} if $F'_{\beta} \approx\varphi$.

\begin{proposition}\label{prop_F}
Let $F\subset C_{p}(X)$ be a vector subspace endowed with the topology induced from $C_{p}(X)$.
If $F$ is $\varphi$-dual, then $F$ is a metrizable and separable space.
\end{proposition}
\begin{proof}
By our assumptions, the strong dual  $F'_{\beta}$ of $F$ is isomorphic to $\varphi$. The strong bidual $F''_{\beta}$ of $F$ is isomorphic to $\mathbb{R}^\omega$,  which does not admit a weaker Hausdorff vector topology, see \cite[Corollary 2.6.5]{bonet}.  Hence the weak*-topology $\sigma(F'',F')$ on $F''$  coincides with the original strong topology of  $F''_{\beta}$.  Consequently,  $(F'',\sigma(F'',F'))$ (so  also $\supset (F,\sigma(F,F'))$) is metrizable and   separable. Since the original topology $\tau_{p}$ of $C_{p}(X)$ is the weak topology, we deduce that  $\tau_{p}|F=\sigma(F,F')$ is  metrizable and separable.
\end{proof}

\begin{theorem}\label{product-1}
 The following assertions are equivalent:
\begin{enumerate}
\item [(1)] $C_{p}(X)$ contains a  $\varphi$-dual  subspace  complemented in $C_{p}(X)$.
\item [(2)] Either $C_{p}([1,\omega])$ or $C_p(\N)$ is a continuous linear image of $C_{p}(X)$.
\end{enumerate}
\end{theorem}
\begin{proof}
(1) $\Rightarrow$  (2): $C_{p}(X)$ contains some complemented subspace $F$ with $F'_{\beta}\approx\varphi$.

\emph{Case 1. $X$ is pseudocompact}, i.e. every $f\in C(X)$ is bounded on $X$. We know that $F$ is a metrizable and separable space,
 by Proposition \ref{prop_F}. In this situation we can apply \cite[Theorem 1]{BKS1} (implication (5) $\Rightarrow$ (4))
and conclude that there exists a continuous linear map from $C_{p}(X)$  onto  $(c_{0})_{p}$, where
$(c_{0})_{p}=\{(x_n)_{n\in\omega}\in \mathbb{R}^\omega\colon x_n\to 0\}$, and   $(c_{0})_{p}\approx C_{p}([1,\omega])$.

\emph{Case 2. $X$ is not pseudocompact}. By  \cite[Section 4]{arkh4},  $C_{p}(X)$ contains a complemented copy of $C_{p}(\N)$,
 so $C_{p}(X)$ can be mapped onto $C_{p}(\N)$.

(2) $\Rightarrow$ (1): Assume that $T: C_{p}(X)\rightarrow C_{p}(\N)$  is a continuous linear surjection.
 If every $g\in C(X)$ is bounded, then $C_{p}(X)=\bigcup_{n\in \omega}\{g\in C(X): |g(x)|\leq n\}$.
 Hence $T(C_{p}(X))$ would be covered by a sequence of bounded sets, which is impossible by the Baire Category Theorem. Therefore $C_{p}(X)$
contains a complemented copy of $\mathbb{R}^\omega \approx C_p(\N)$, by \cite[Section 4]{arkh4}.
Finally, assume that $C_{p}([1,\omega])$ is a continuous linear image of $C_{p}(X)$. By \cite[Theorem 1]{BKS1} (implication (4) $\Rightarrow$ (2)),
$C_{p}(X)$ contains a complemented copy of $(c_{0})_{p}\approx C_{p}([1,\omega])$.
In all cases $C_p(\N)$ or $C_{p}([1,\omega])$ are $\varphi$-dual subspaces of $C_{p}(X)$.
\end{proof}

\begin{corollary}\label{XxY}
For infinite Tychonoff spaces $X$ and $Y$ the lcs
 $C_{p}(X\times Y)$ always contains a complemented $\varphi$-dual subspace.
\end{corollary}
\begin{proof}
 If $X\times Y$ is \emph{not pseudocompact}, then $C_{p}(X\times Y)$ contains a complemented copy of $\mathbb{R}^{\omega} \approx C_p(\N)$.
 If $X\times Y$ is \emph{pseudocompact}, then $C_{p}(X\times Y)$ contains a complemented copy of  $(c_{0})_{p}\approx C_{p}([1,\omega])$, by \cite[Theorem 1.4]{KSMZ}.
\end{proof}

In view of Theorem \ref{some-remark} the following questions arise naturally.

\begin{problem}
Does the implication (1) $\rightarrow$ (2)  in Theorem \ref{some-remark} remain valid
without an assumption that $C_{p}(X)$ is distinguished?
\end{problem}

\begin{problem} Does the implication (3) $\rightarrow$ (2) in Theorem \ref{some-remark} remain valid
without an assumption that $C_{p}(X)$ is distinguished? More specifically,
does there exist a non-open continuous linear surjection $T:C_p[0,1]\rightarrow C_p[0,1]$ satisfying condition (3) of Theorem \ref{some-remark}?
\end{problem}
\section{Distinguished spaces $C_p([1,\alpha])$ for countable ordinals $\alpha$}\label{sec3}

In this section we are interested in finding a complete characterization of those spaces $Y$ which admit a continuous linear surjection $T:C_{p}(X)\rightarrow C_p(Y)$
for some simplest $\Delta$-spaces $X$.

1) Let $X$ be a discrete space. It follows from our \cite[Corollary 3.3]{Kakol-Leiderman} that a continuous linear surjection $T:C_{p}(X)\rightarrow C_p(Y)$
exists if and only if $Y$ itself is discrete and $|Y| \leq |X|$.

2) Assume that $X$ is  a metrizable compact space. Then, as we showed in \cite[Proposition 3.5]{ka-arkady},
$X$ is a $\Delta$-space if and only if $X$ is countable, equivalently, if and only if $X$ is homeomorphic to the ordinal space $[1,\alpha]$ for some countable ordinal $\alpha$,
by the Mazurkiewicz-Sierpi\'nski theorem.
If we drop the requirement of linearity of the operator $T$, then for every metrizable compact space $Y$
there exists a continuous surjection $T: C_{p}([1,\omega]) \rightarrow C_{p}(Y)$, see \cite[Proposition 5.4]{Krupski} and \cite[Remark 3.4]{Kawamura}.

Let us fix an infinite countable ordinal $\alpha$ and ask for which ordinal $\beta  \geq \alpha$ there exist a continuous linear isomorphism
$T: C_{p}([1,\alpha]) \rightarrow C_{p}([1,\beta])$?
The answer for this question has been given independently in \cite{Baars} and \cite{Gul'ko}.

\begin{theorem} \cite{Baars}, \cite{Gul'ko} \label{th:BG}
Let $\omega \leq \alpha \leq \beta < \omega_1$. Then the following conditions are equivalent:
\begin{enumerate}
\item [(1)] $C_{p}([1,\alpha])$ and $C_{p}([1,\beta])$ are isomorphic.
\item [(2)] The Banach spaces $C([1,\alpha])$ and $C([1,\beta])$ are isomorphic.
\item [(3)] $\beta < \alpha^{\omega}$.
\end{enumerate}
\end{theorem}

Note that the most difficult implication (2) $\rightarrow$ (3) in the proof of Theorem \ref{th:BG} in both papers \cite{Baars} and \cite{Gul'ko}
is based on the result of C. Bessaga and A. Pe{\l}czy\'nski \cite[Lemma 2]{bessaga}.
However, the proof of this key Lemma 2 presented in \cite{bessaga} heavily depends on the assumption that there exists a
linear isomorphism $T: C([1,\alpha]) \to C([1,\beta])$ and it does not work if we assume only that $T$ is a continuous linear surjection even for $\alpha=\omega$.
This obstacle has been resolved in the following new result.

\begin{theorem}\label{cpp}
Let $\alpha$ be a fixed infinite countable ordinal. Then for a Tychonoff space $Y$ the following are equivalent.
\begin{enumerate}
\item [(1)] There exists a linear continuous surjection $T:C_p([1,\alpha])\rightarrow C_p(Y)$.
\item [(2)] $Y$ is homeomorphic to $[1,\beta]$, where $\beta$ is a countable ordinal such that either $\beta < \alpha$, or $\alpha \leq \beta < \alpha^{\omega}$.
\item [(3)] $Y$ is homeomorphic to $[1,\beta]$, where $\beta$ is a countable ordinal with the property: if $\alpha \leq \beta$ then
$C_{p}([1,\alpha])$ and $C_{p}([1,\beta])$ are isomorphic.
 \end{enumerate}
 \end{theorem}

\begin{proof} Implication (2) $\Rightarrow$ (3) is an immediate consequence of Theorem \ref{th:BG}. If $\beta < \alpha$ then we define the operator of restriction:
$T:f\mapsto f|_{[1,\beta]}$. If $\alpha \leq \beta$, then the assumption of (3) says that there is an isomorphism between $C_{p}([1,\alpha])$ and $C_{p}([1,\beta])$
which trivially implies (3) $\Rightarrow$ (1). Now we prove (1) $\Rightarrow$ (2).

It is probably a folklore that (1) implies that $Y$ has to be a countable compact space. For instance, one can argue as follows.
First, $Y$ is a metrizable compact space by \cite[Theorem 2.7]{Leiderman}. Second, $Y$ is scattered by \cite[Proposition 3.9]{Kakol-Leiderman}.
Every scattered metrizable compact space is homeomorphic to $[1,\beta]$ for some countable ordinal $\beta$. If $\beta < \alpha$ there is nothing to prove.
So, let us assume that $\alpha \leq \beta$.

Applying the Closed Graph Theorem we consider $T$ as a linear bounded operator from the Banach space $C([1,\alpha])$ onto the Banach space $C([1,\beta])$.
Note that in a particular case $\alpha= \omega$ we have that $C([1,\alpha]) \approx c_{0}$, hence $C([1,\beta])$ is a quotient of $c_0$,
and we can apply the theorem of W.B. Johnson and M. Zippin \cite{Johnson} to deduce that the space $C([1,\beta])$ is isomorphic to $c_0$. By Theorem \ref{th:BG} we conclude that in
this particular case $\beta < \omega^{\omega}$.

In a general case of arbitrary countable ordinal $\alpha$ we will apply more elaborate results of functional analysis.
Recall that the Szlenk index of a Banach space $E$, denoted ${\rm Sz}(E)$, is an ordinal number, which is invariant under linear isomorphisms.
For the definition of ${\rm Sz}(E)$ and its basic properties we refer the reader to the survey papers \cite{Lancien}, \cite{Rosenthal}.
The key tool is the following precise result of C. Samuel, which is in turn based on a deep result of D. Alspach and Y. Benyamini \cite{Alspach},
(see \cite[Theorem 7]{Lancien}, \cite[Theorem 2.15]{Rosenthal}).

\noindent
\textbf{Fact A.} \textit{For any} $0 \leq \gamma <\omega_1$
$${\rm Sz}(C([1,\omega^{\omega^{\gamma}}])) = \omega^{\gamma+1}.$$

We need also

\noindent
\textbf{Fact B.} \cite[Corollary 2.19]{Rosenthal} \textit{Let $E_1$ and $E_2$ be given Banach spaces with norm-separable duals.
Assume that $E_2$ is isomorphic to a subspace of a quotient space of $E_1$.
Then ${\rm Sz}(E_2) \leq {\rm Sz}(E_1)$.}

\medskip
In order to finish the proof of (1) $\Rightarrow$ (2) suppose the contrary: $\beta \geq \alpha^{\omega}$.
 Then by Fact A, ${\rm Sz}(C([1,\beta])) > {\rm Sz}(C([1,\alpha]))$ which contradicts Fact B.
\end{proof}

\begin{remark}
Linear continuous surjections $T:C_p([1,\alpha])\rightarrow C_p([1,\beta])$ constitute only a proper part of the set of all bounded operators between corresponding
Banach spaces. It would be very helpful to find a more direct topological argument for the proof of $\beta < \alpha^{\omega}$ in (1) $\Rightarrow$ (2) above.
\end{remark}

An lcs $E$ is called \emph{hereditarily distinguished} if every closed linear subspace of $E$ is distinguished \cite{Kakol-Leiderman}.
Every Tychonoff product $\R^X$ is hereditarily distinguished.
However, it is known that even a Fr\'echet distinguished lcs can contain a closed non-distinguished subspace.
Problem 3.18 \cite{Kakol-Leiderman} asks: Does there exist an infinite compact space $X$ such that $C_p(X)$ is hereditarily distinguished?
Below we answer in the affirmative to this question.

\begin{theorem}\label{th:hered}
The space $C_p([1,\alpha])$ is hereditarily distinguished for every countable ordinal $\alpha$.
\end{theorem}
\begin{proof}
Denote $X= [1,\alpha]$. Then the product $\R^X$ is metrizable, since $X$ is countable.
Let $E$ be any closed linear infinite-dimensional subspace of $C_p(X)$. If $f \in \R^X$ belongs to the closure of $E$ in $\R^X$, then there is a sequence $(f_n)$ in $E$
converging to $f$. This sequence provides a bounded set in $E$. This means that $E$ is a dense large subspace of $\R^Y$ for some infinite countable set $Y$.
Then the strong duals of both spaces $E$ and $\R^Y$ are the same and carry the finest locally convex topology since the strong dual of $\R^Y$ is $\phi$.
\end{proof}

\section{From pointwise topology to weak topology}\label{sec4}

M. Krupski proved that if $X$ is an infinite metrizable $C$-space then $C_p(X)$ and $C_{k}(X)_{w}$ are not homeomorphic \cite{Krupski-1}.
Recall that a normal space $X$ is said to be a  \emph{$C$-space} if for any given sequence of its open covers $\{\ucal_n:n\in\omega\}$
there exists a sequence $\{\vcal_n:n\in\omega\}$ of families of pairwise disjoint open sets such that
$\vcal_n$ is a refinement of $\ucal_n$ for each $n\in\omega$, and $\bigcup_{n\in\omega} \vcal_n$ is a cover of $X$.

Every countable-dimensional (in particular, every finite-dimensional space) is a $C$-space.
It is well known that the Hilbert cube $Q=[0,1]^\omega$ is not a $C$-space.

Furthermore,  M. Krupski and W. Marciszewski proved that if $K$ and $L$ are infinite compact spaces, then there is no a homeomorphism
$\Phi: C_{k}(K)_{w} \rightarrow C_{p}(L)$, which is in addition uniform (see \cite[Proposition 3.1]{Krupski-2}). In particular,
a linear homeomorphism between $C_{p}(L)$ and $C_{k}(K)_{w}$ does not exist.
Note that the proof presented in \cite{Krupski-2} uses very essentially the assumption that $\Phi^{-1}(A)$ is a compact subset of $C_{k}(Y)_{w}$ provided
 $A$ is a converging sequence in $C_{p}(L)$, and it is not clear whether this assumption can be removed.
These remarks explain our motivation for the next open question.

\begin{problem}\label{pro-1}
Does there exist an infinite compact space $X$ admitting a continuous linear surjection $T:C_{p}(X)\rightarrow C_{k}(X)_{w}$?
\end{problem}

In Theorem \ref{the-1} we show that such $X$ (if exists) cannot be an infinite metrizable compact $C$-space.
In particular, $X$ cannot be an infinite metrizable finite-dimensional compact space, equivalently, an infinite compact subspace of the Euclidean space $\R^n$.

\begin{observation}\label{obs}
If $Y$ is an infinite Tychonoff space and $T:C_{p}(X)\rightarrow C_{k}(Y)_{w}$ is a continuous open linear surjection,
 then the topologies $\tau_{k}$ and $w$ coincide on $C(Y)$.
Hence such an operator $T$ does not exist if $Y$ is an infinite compact space.
\end{observation}
 \begin{proof} Indeed, let  $\widehat{T}: C_{p}(X)/\ker(T) \rightarrow C_{k}(Y)_{w}$ be an isomorphism  associated with $T$, where $\widehat{T}\circ P=T$ and
 $P:C_{p}(X)\rightarrow C_{p}(X)/\ker(T)$ is the quotient map.
Since $C_{p}(X)$ is quasibarrelled (hence also $C_{p}(X)/\ker(T)$), the space $C_{k}(Y)_{w}$ is quasibarrelled. But the topology $\tau_{k}$ is $w$-polar, so $w=\tau_{k}$.
\end{proof}

Assuming that $X$ in the above Claim \ref{obs} is a $\Delta$-space one can remove a restriction that $T$ is open.

\begin{theorem}\label{Kru-Mar}
Let $T:C_{p}(X)\rightarrow C_{k}(Y)_{w}$ be a continuous linear surjection, where $X$ is a $\Delta$-space and $Y$ is a Tychonoff space.
Then  $C_{k}(Y)_{w}=C_{p}(Y)=C_{k}(Y)$. Consequently, if $Y$ is an infinite compact space and $X$ is a $\Delta$-space, then such an operator $T$ does not exist.
\end{theorem}

\begin{proof}
Following  the first part of the proof of (1) $\Rightarrow$ (2) in Theorem \ref{Ck} the strong dual of $C_{k}(Y)_{\sigma}$ is free.
 Since $C_{k}(Y)$ and $C_{k}(Y)_{w}$ have the same bounded sets and the same topological dual, $C_{k}(Y)'_{\beta}$ is free.
 By Theorem \ref{Ck}, we obtain that $\tau_{p}=\tau_{k}$ on $C(Y)$. This is impossible for infinite compact space $Y$.
\end{proof}

\begin{theorem}\label{the-1}
For an infinite metrizable compact $C$-space $X$ a continuous linear surjection $T:C_{p}(X)\rightarrow C_{k}(X)_{w}$ does not exist.
\end{theorem}
\begin{proof}
On the contrary, assume that such $X$ and $T$ exist. First we eliminate a possibility that $X$ is countable.
Indeed, every countable $X$ is a $\Delta$-space and we can apply Theorem \ref{Kru-Mar}.
Let $X$ be uncountable, then by the celebrated Milyutin theorem the Banach spaces $C_{k}(Q)$ and $C_{k}(X)$ are isomorphic,
where $Q$ denotes the Hilbert cube,
 so also $C_{k}(Q)_{w}$ and $C_{k}(X)_{w}$ are isomorphic.
Since the identity mapping from $C_{k}(Q)_{w}$ onto $C_p(Q)$ is continuous, we get a continuous linear surjection
$S: C_p(X)\rightarrow C_p(Q)$.
Our next argument is based on an analysis of the dual spaces $C_p(X)'=L_p(X)$
(for more details see proofs of \cite[Proposition 2.1]{Kawamura} and \cite[Theorem 3.4]{Leiderman-1}).

The dual mapping $S^*$ embeds $L_p(Q)$ into $L_p(X)$, and therefore
 $Q$ is homeomorphic to a subspace of $L_p(X)$.
For each natural $n\in\N$ consider the subspace $B_n(X)$
of $L_p(X)$ formed by all words of the reduced
length $\leq n$ over $X$, and let $A_n(X) = B_n(X) \setminus B_{n-1}(X)$ be the subspace of all words of the reduced length $n$.
It is known that each $B_n(X)$ is closed in $L_p(X)$, and $A_n(X)$ is homeomorphic to the subspace
$(\R^*)^n \times (X^n \setminus \Lambda_n)$ of the Cartesian product $(\R^*)^n \times X^n$,  where $\R^* = \R\setminus\{0\}$
and $\Lambda_n = \{(x_1, \dots, x_n)\in X^n: x_i = x_j$ for some $i \neq j\}$.
Obviously,  $\Lambda_n$ is closed in the metrizable space $X^n$,
hence $\Lambda_n$ is a $G_{\delta}$-set in $X^n$ and $A_n(X)$ is a $F_{\sigma}$-set in $(\R^*)^n \times X^n$.
Each Cartesian product $(\R^*)^n \times X^n$ is a $C$-space by \cite[Theorem 3]{Rohm},
then each $A_n(X)$ is a $C$-space by \cite[Subspace Theorem]{Rohm}, hence
$L_p(X) = \bigcup_{n\in\omega} A_n(X)$ is a $C$-space by \cite[Sum Theorem]{Rohm}
and finally, the Hilbert cube $Q$ would be also a $C$-space, which is false.
The obtained contradiction finishes the proof of Theorem \ref{the-1}.
\end{proof}

\begin{problem}\label{product}
Find a characterization of metrizable compact spaces $X$ admitting  a continuous linear surjection $T: C_p(X)\rightarrow C_p(Q)$.
\end{problem}

\begin{remark}\label{Daverman}
R. Daverman kindly informed us that apparently there should exist a metrizable compact space $X$ with the following properties:
the square $X^2$ is homeomorphic to the Hilbert cube $Q$, but $X$ itself does not contain a homeomorphic copy of $Q$.
By these reasons a simple answer to Problem \ref{product} seems to be unavailable. Obviously, $X$ cannot be a $C$-space, by the proof of Theorem \ref{the-1}.
Note also that a continuous linear surjection $T: C_p[0,1]\rightarrow C_p(X)$ exists for any finite-dimensional metrizable compact space $X$ \cite{Leiderman}.
\end{remark}

\end{document}